\documentclass[11pt]{amsart}
\usepackage[latin1]{inputenc}
\usepackage[T1]{fontenc}
\usepackage{hyperref}
\usepackage{amsfonts}
\usepackage{amsmath}
\usepackage{amstext}
\usepackage[dvips]{graphicx}
\usepackage{epsf, subfigure, verbatim}
\usepackage{stmaryrd}
\usepackage{amssymb}
\usepackage{amsthm}
\usepackage{latexsym}
\usepackage{layout}
\usepackage{xcolor}

\newtheorem{thm}{Theorem}[section]
\newtheorem{prop}[thm]{Proposition}
\newtheorem{lem}[thm]{Lemma}

\newtheorem{assump}[thm]{Assumption}

\newtheorem{rem}[thm]{Remark}
\newtheorem{defn}[thm]{Definition}

\newcommand{\reels}{\mathbb{R}}
\newcommand{\nat}{\mathbb{N}}
\newcommand{\ind}{\mathbb{I}}
\newcommand{\esp}{\mathbb{E}}
\newcommand{\p}{\mathbb{P}}
\newcommand{\tribu}{\mathcal{F}}
\newcommand{\tribuborel}{\mathcal{B}(\reels)}
\newcommand{\tribuprev}{\mathcal{P}}
\newcommand{\filt}{\mathcal{F}_t}
\newcommand{\D}{\mathbb{D}^{1,2}}
\newcommand{\mesleb}{\lambda}
\newcommand{\meslevy}{\nu}
\newcommand{\intPoissont}{\int_0^t \int_{\reels}}
\newcommand{\intPoissontmoins}{\int_0^{t-} \int_{\reels}}
\newcommand{\intPoissonT}{\int_0^T \int_{\reels}}
\newcommand{\Ndt}{N(dt,dz)}
\newcommand{\Nds}{N(ds,dz)}
\newcommand{\Ntildedt}{\widetilde{N}(dt,dz)}
\newcommand{\Ntildeds}{\widetilde{N}(ds,dz)}

\newcommand{\Ntildedyds}{\widetilde{N}(ds,dy)}
\newcommand{\Ntildedydr}{\widetilde{N}(dr,dy)}
\newcommand{\x}{\mathcal{X}}
\newcommand{\sigman}{\Sigma_{(i_1,\dots,i_n)} (\x)}

\newcommand{\sigmam}{\Sigma_{(j_1,\dots,j_m)} (\x)}

\newcommand{\sigmanomisk}{\Sigma_{(i_1,\dots,\widehat{i_k},\dots,i_n)} (\x)}
\newcommand{\zerotn}{\Pi_{i_1}(\x) \times \dots \times \Pi_{i_n}(\x)}

\begin{document}

\title[Malliavin calculus for Lévy processes]{Malliavin calculus and Clark-Ocone formula for functionals of a square-integrable Lévy process}

\author{Jean-François Renaud}
\thanks{Corresponding author: renaud@dms.umontreal.ca}
\address{J.-F. Renaud: D{\'e}partement de math{\'e}matiques et de statistique, Universit{\'e} de Montr{\'e}al,
C.P. 6128, Succ. Centre-Ville, Montr{\'e}al, Qu{\'e}bec, H3C 3J7, Canada}
\email{renaud@dms.umontreal.ca}

\author{Bruno R{\'e}millard}
\address{B. Rémillard: Service de l'enseignement des m{\'e}thodes quantitatives de gestion, HEC Montr{\'e}al,
3000 chemin de la C{\^o}te-Sainte-Catherine, Montr{\'e}al, Qu{\'e}bec, H3T 2A7, Canada}
\email{bruno.remillard@hec.ca}


\maketitle

\begin{abstract}
In this paper, we construct a Malliavin derivative for functionals of square-integrable Lévy processes and derive a Clark-Ocone formula. The Malliavin derivative is defined via chaos expansions involving stochastic integrals with respect to Brownian motion and Poisson random measure. As an illustration, we compute the explicit martingale representation for the maximum of a Lévy process.
\end{abstract}

\section{Introduction}

If $W = (W_t)_{t \in [0,T]}$ is a Brownian motion, then the Wiener-Itô chaos expansion of a square-integrable Brownian functional $F$ is given by
\begin{equation}\label{E:wienerito}
F = \esp[F] + \sum_{n \geq 1} \int_0^T \dots \int_0^T f_n(t_1,\dots,t_n) \, W(dt_1) \dots W(dt_n) ,
\end{equation}
where $(f_n)_{n \geq 1}$ is a sequence of deterministic functions. This chaotic representation can be obtained by iterating Itô's representation theorem and can then be used to define the classical Malliavin derivative in the following way: if the chaos expansion of $F$ satisfies an integrability condition, then $F$ is Malliavin-differentiable and its Malliavin derivative $DF$ is given by
\begin{multline}\label{E:derMalliavin}
D_t F = f_1(t) \\
+ \sum_{n \geq 1} (n+1) \int_0^T \dots \int_0^T f_{n+1} (t_1,\dots,t_n,t) \, W(dt_1) \dots W(dt_n) ,
\end{multline}
for $t \in [0,T]$. This derivative operator is equal to a weak derivative on the Wiener space; the close connection between Hermite polynomials and Brownian motion is at the hearth of that equivalence. See for instance Nualart \cite{nualart1995}.

Quite recently, L{\o}kka \cite{lokka2004} developed similar results for a square-integrable pure-jump Lévy process $L = (L_t)_{t \in [0,T]}$ given by
$$
L_t = \int_0^t \int_{\reels} z (\mu - \pi)(ds,dz) ,
$$
where $\mu - \pi$ is the compensated Poisson random measure associated with $L$. In this setup, by mimicking the steps of the Wiener-Itô expansion, L{\o}kka obtained a chaos representation property for the pure-jump Lévy process $L$ just as in Equation~\eqref{E:wienerito} and then defined the corresponding Malliavin derivative as in Equation~\eqref{E:derMalliavin}. Later on, Benth et al. \cite{benthetal2003} introduced chaos expansions and a Malliavin derivative for more general Lévy processes, i.e. Lévy processes with a Brownian component. However, in the latter, neither proofs nor connections with the classical definitions was given.

Our first goal is to provide a detailed construction of a chaotic Malliavin derivative leading to a Clark-Ocone formula for Lévy processes. We extend the definitions of the Malliavin derivatives for Brownian motion and pure-jump Lévy processes to general square-integrable Lévy processes. Secondly, we derive additional results that are useful for computational purposes.

Our approach follows more or less the same steps as those leading to the Wiener-Itô chaos expansion and the chaotic Brownian Malliavin derivative, just as L{\o}kka \cite{lokka2004} did for pure-jump Lévy processes. The definition of the directional Malliavin derivatives is different from those of Benth et al. \cite{benthetal2003}. The main idea is to obtain a chaotic representation property (CRP) by iterating a well-chosen martingale representation property (MRP) and then defining directional Malliavin derivatives in the spirit of Ma et al. \cite{maetal1998}. However, in the context of a general square-integrable Lévy process, one has to deal with two integrators and therefore must be careful with the choice of derivative operators in order to extend the classical definitions. This choice will be made with the so-called commutativity relationships in mind and following Le{\'o}n et al. \cite{leonetal2002}. In the Brownian motion setup, the commutativity relationship between Malliavin derivative and Skorohod integral is given by
\begin{equation}\label{E:commute}
D_t \int_0^T u_s \, W(ds) = u_t + \int_t^T D_t u_s \, W(ds) ,
\end{equation}
when $u$ is an adapted process. See Theorem 4.2 in Nualart and Vives \cite{nualartvives1990} for the corresponding formula in the Poisson process setup.

We will get the MRP using a denseness argument involving Doléans-Dade exponentials. Our path toward the CRP is different from that of Itô \cite{ito1956} and Kunita and Watanabe \cite{kunitawatanabe1967} who used random measures; see also the recent formulation of that approach given by Kunita \cite{kunita2004} and Solé et al. \cite{soleetal2007}. It is known that the CRP usually implies the MRP and that in general a Lévy process does not possess the MRP nor a predictable representation property. However, we show that the CRP and our well-chosen MRP are equivalent for square-integrable Lévy processes. Finally, just as in the Brownian and pure-jump Lévy setups, a Malliavin derivative and a Clark-Ocone formula are derived. As an illustration, we compute the explicit martingale representation for the maximum of a Lévy process.

This approach to Malliavin calculus for Lévy processes is different from the very interesting contributions of Nualart and Schoutens \cite{nualartschoutens2000}, Le{\'o}n et al. \cite{leonetal2002} and Davis and Johansson \cite{davisjohansson2006}. They developed in sequence a Malliavin calculus for Lévy processes using different chaotic decompositions based on orthogonal polynomials. Their construction also relies on the fact that all the moments of their Lévy process exist. Many other chaos decompositions related to Lévy processes have been considered through the years: see for example the papers of Dermoune \cite{dermoune1990}, Nualart and Vives \cite{nualartvives1990}, Aase et al. \cite{aaseetal2000} and Lytvynov \cite{lytvynov2003}.

On the other hand, Kulik \cite{kulik2006} developed a Malliavin calculus for Lévy processes in order to study the absolute continuity of solutions of stochastic differential equations with jumps, while Bally et al. \cite{ballyetal2007} established an integration by parts formula in order to give numerical algorithms for sensitivity computations in a model driven by a Lévy process; see also Bavouzet-Morel and Messaoud \cite{bavouzetmessaoud2006}. Finally, in a very interesting paper, Solé et al. \cite{soleetal2007} constructed a Malliavin calculus for Lévy processes through a suitable canonical space. While finishing this paper, the work of Petrou \cite{petrou2006} was brought to our attention. In that paper, the same methodology is applied to obtain a Malliavin derivative and a Clark-Ocone formula, but the focus is on financial applications. Since one of our goal is to give a thorough treatment of a Malliavin calculus for square-integrable Lévy processes, and considering the major differences between the two papers, we think that each paper has his own interest.

The rest of the paper is organized as follows. In Section 2, preliminary results on Lévy processes are recalled. In Section 3 and 4, martingale and chaotic representations are successively obtained. Then, in Section 5, the corresponding Malliavin derivative is constructed in order to get a Clark-Ocone formula. Finally, in Section 6, we apply this Clark-Ocone formula to compute the martingale representation of the maximum of a Lévy process.

\section{Preliminary results on Lévy processes}

Let $T$ be a strictly positive real number and let $X = (X_t)_{t \in [0,T]}$ be a Lévy process defined on a probability space $(\Omega, \tribu, \p)$, i.e. $X$ is a process with independent and stationary increments, is continuous in probability and starts from $0$ almost surely. We assume that $X$ is the \textit{càdlàg} modification and that the probability space is equipped with the completed filtration $(\filt)_{t \in [0,T]}$ generated by $X$. We also assume that the $\sigma$-field $\tribu$ is equal to $\tribu_T$.

This filtration satisfies \textit{les conditions habituelles} and, for any fixed time $t$, $\tribu_{t-} = \tribu_t$. Consequently, the filtration is continuous. This fact is crucial in the statement of our Clark-Ocone formula.

The reader not familiar with Lévy processes is invited to have a look at the books of Schoutens \cite{schoutens2003}, Protter \cite{protter2004} and Bertoin \cite{bertoin1996}.

From the Lévy-Itô decomposition (see \cite{protter2004}, Theorem $42$), we know that $X$ can be expressed as
\begin{equation}\label{E:formegen}
X_t = \alpha t + \sigma W_t + \int_0^t \int_{|z| \geq 1} z \, \Nds + \int_0^t \int_{|z| < 1} z \, \Ntildeds
\end{equation}
where $\alpha$ is a real number, $\sigma$ is a strictly positive real number, $W$ is a standard Brownian motion and $\widetilde{N}$ is the compensated Poisson random measure associated with the Poisson random measure $N$. The Poisson random measure $N$ is independent of the Brownian motion $W$. Its compensator measure is denoted by $\mesleb \times \meslevy$, where $\mesleb$ is Lebesgue measure on $[0,T]$ and $\meslevy$ is the Lévy measure of $X$, i.e. $\meslevy$ is a $\sigma$-finite measure on $\reels$ such that $\meslevy(\{0\}) = 0$ and
$$
\int_{\reels} (1 \wedge z^2) \, \meslevy(dz) < \infty .
$$
Therefore the compensated random measure $\widetilde{N}$ is defined by
$$
\widetilde{N} ([0,t] \times A) = N ([0,t] \times A) - t \meslevy(A) .
$$
This measure is equal to the measure $\mu - \pi$ mentioned in the introduction.

Finally, let $\tribuprev$ be the predictable $\sigma$-field on $[0,T] \times \Omega$ and $\tribuborel$ the Borel $\sigma$-field on $\reels$. We recall that a process $\psi(t,z,\omega)$ is Borel predictable if it is $(\tribuprev \times \tribuborel)$-measurable.

\subsection{Square-integrable Lévy processes}

When the Lévy process $X$ is square-integrable, it can also be expressed as
\begin{equation}\label{E:Lgen}
X_t = \mu t + \sigma W_t + \intPoissont z \, \Ntildeds ,
\end{equation}
where $\mu = \esp [X_1]$. Indeed, in Equation~\eqref{E:formegen} we have that
$$
\alpha = \esp \left[X_1 - \int_0^1 \int_{|z| \geq 1} z \, \Ndt \right] ,
$$
so $\esp[X_t^2]$ is finite if and only if
$$
\int_{\reels} z^2 \meslevy(dz) = \esp \left[ \left( \int_0^1 \int_{|z| \geq 1} z \, \Ndt \right)^2 \right]
$$
is finite. Note that in general $\mu \neq \alpha$.

Here is a consequence of Itô's formula.
\begin{lem}\label{L:procZ}
If $h$ belongs to $L^2([0,T], \mesleb)$ and if $(t,z) \mapsto e^{g(t,z)} - 1$ belongs to $L^1([0,T] \times \reels, \mesleb \times \meslevy)$, define $Z = (Z_t)_{t \in [0,T]}$ by
\begin{multline}\label{E:procZ}
Z_t = \exp \left\lbrace \int_0^t h(s) \, W(ds) - \frac{1}{2} \int_0^t h^2(s) \, ds + \intPoissont g(s,z) \, \Nds \right. \\
- \left. \intPoissont \left( e^{g(s,z)} - 1 \right) \, \meslevy(dz) ds \right\rbrace .
\end{multline}
The process $Z$ is a square-integrable martingale if and only if $e^g - 1$ is an element of $L^2([0,T] \times \reels, \mesleb \times \meslevy)$.
\end{lem}
\begin{proof}
From the assumptions, we have that $g$ belongs to $L^2([0,T] \times \reels, \mesleb \times \meslevy)$ and that $Z$ is a well-defined positive local martingale. Then, if $\esp[Z_T] = 1$, it is a martingale. From Itô's formula, we also have that $Z$ is the solution of
$$
dZ_t = Z_{t-} \, h(t) \, W(dt) + Z_{t-} \int_{\reels} (e^{g(t,z)} - 1) \, \Ntildedt , \quad Z_0 = 1 .
$$
Let $(\tau_n)_{n \geq 1}$ be the fundamental sequence of stopping times of $Z$. Since $W$ and $N$ are independent,
$$
\esp [Z_{t \wedge \tau_n}^2] = 1 + \esp \left[ \int_0^{t \wedge \tau_n} Z_s^2 \, h^2(s) \, ds \right] + \esp\left[ \int_0^{t \wedge \tau_n} Z_s^2 \int_{\reels} (e^{g(s,z)} - 1)^2 \, \meslevy(dz) ds \right] ,
$$
for every $n \geq 1$. Taking the limit when $n$ goes to infinity yields
\begin{equation}\label{E:egalitequadratique}
\esp [Z_t^2] = 1 + \int_0^t \esp[Z_s^2] \, h^2(s) \, ds + \int_0^t \esp[Z_s^2] \int_{\reels} (e^{g(s,z)} - 1)^2 \, \meslevy(dz) ds .
\end{equation}
If we define $G(t) = h^2(t) + \int_{\reels} (e^{g(t,z)} - 1)^2 \, \meslevy(dz)$, then the function $t \mapsto \esp [Z_t^2]$ is the solution of
$$
F^{\prime}(t) = G(t) F(t) , \quad F(0) = 1 .
$$
Hence,
\begin{equation}\label{E:lemprocZ}
\esp [Z_t^2] = \exp \left\lbrace \int_0^t h^2(s) \, ds + \intPoissont (e^{g(s,z)} - 1)^2 \, \meslevy(dz) ds \right\rbrace
\end{equation}
and the statement follows.
\end{proof}

For $h \in L^2([0,T], \mesleb)$ and $e^g - 1 \in L^2([0,T] \times \reels, \mesleb \times \meslevy)$, the process $Z$ is the Doléans-Dade exponential of the square-integrable martingale $(\overline{M}_t)_{t \in [0,T]}$ defined by
$$
\overline{M}_t = \int_0^t h(s) \, W(ds) + \int_0^t \int_{\reels} (e^{g(s,z)} - 1) \, \Ntildeds .
$$
In the literature, this is often denoted by $Z = \mathcal{E}(\overline{M})$, the stochastic exponential of $\overline{M}$.

\subsection{A particular choice for $g$}

If $g$ is an element of $L^2([0,T] \times \reels, \mesleb \times \meslevy)$, then $e^{g} - 1$ is not necessarily square-integrable. One way to circumvent this problem is to introduce the bijection $\gamma \colon \reels \to (-1,1)$ defined by
\begin{equation}\label{E:gamma}
\gamma(z) =
\begin{cases}
e^z - 1 \quad \text{if $z < 0$,}\\
1 - e^{-z} \quad \text{if $z \geq 0$.}
\end{cases}
\end{equation}
Note that $\gamma$ is bounded. Hence, if $h$ is square-integrable on $[0,T]$ and if $g$ is of the form $g(t,z) = \bar{g}(t) \gamma(z)$, where $\bar{g} \in C([0,T])$, i.e. $\bar{g}$ is a continuous function on $[0,T]$, then $Z$ is square-integrable by Lemma~\ref{L:procZ}.

The idea of introducing the function $\gamma$ is taken from L{\o}kka \cite{lokka2004}. In that paper, it is also proved that the process $(N_t)_{t \in [0,T]}$ defined by
\begin{equation}\label{E:defN}
N_t = \intPoissont z \, \Ntildeds
\end{equation}
and the process $(\widehat{N}_t)_{t \in [0,T]}$ defined by
$$
\widehat{N}_t = \intPoissont \gamma(z) \, \Ntildeds
$$
generate the same filtration. Since
$$
\tribu_t^X = \tribu_t^W \vee \tribu_t^N
$$
for every $t \in [0,T]$ (see Lemma $3.1$ in \cite{soleetal2007}), we have the following lemma.
\begin{lem}
For every $t \in [0,T]$,
$$
\tribu_t^X = \tribu_t^W \vee \tribu_t^N = \tribu_t^W \vee \tribu_t^{\widehat{N}} .
$$
As a consequence, $\tribu = \tribu_T^W \vee \tribu_T^{\widehat{N}}$.
\end{lem}

This means that the processes $X_t = \mu t + \sigma W_t + N_t$ and $\widehat{X}_t = \mu t + \sigma W_t + \widehat{N}_t$ both generate the filtration $(\filt)_{t \in [0,T]}$.

\section{Martingale representations}

\begin{assump}
For the rest of the paper, we suppose that $X$ is a square-integrable Lévy process with a decomposition as in Equation~\eqref{E:Lgen}.
\end{assump}

In general, a Lévy process does not possess the classical \textit{predictable representation property} (PRP), i.e. an integrable random variable $F$ (even with finite higher moments) can not always be expressed as
\begin{equation*}
F = \esp [F] + \int_0^T u_t \, dX_t ,
\end{equation*}
where $u$ is a predictable process and where the stochastic integral is understood as an integral with respect to a semimartingale. However, a martingale representation property exists for square-integrable functionals of $X$. It is a representation with respect to $W(dt)$ and $\Ntildedt$ simultaneously. This result can be found as far back as the paper of Itô \cite{ito1956}. In this section, we will provide a different proof. But first, here is a preparatory lemma.

\begin{lem}\label{L:densite}
The linear subspace of $L^2(\Omega, \tribu, \p)$ generated by
$$
\left\lbrace Y(h,g) \mid h \in L^2([0,T], \mesleb) , \, g \in C([0,T]) \right\rbrace ,
$$
where the random variables $Y(h,g)$ are defined by
\begin{equation}\label{E:dense}
Y(h,g) = \exp \left\lbrace \int_0^T h(t) \, W(dt) + \intPoissonT g(t) \gamma(z) \, \Ntildedt \right\rbrace ,
\end{equation}
is dense.
\end{lem}
\begin{proof}
We adapt the proof of Lemma $1.1.2$ in the book of Nualart \cite{nualart1995}. Let $X$ be a square-integrable random variable such that
\begin{equation*}
\esp \left[ X Y(h,g) \right] = 0
\end{equation*}
for every $h \in L^2([0,T], \mesleb)$ and $g \in C([0,T])$. Let $W(h) = \int_0^T h(t) \, W(dt)$ and $\widetilde{N}(g) = \intPoissonT g(t) \gamma(z) \, \Ntildedt$. Hence,
\begin{equation*}
\esp \left[ X \exp \left\lbrace \sum_{i=1}^n \left( a_i W(h_i) + b_i \widetilde{N}(g_i) \right) \right\rbrace \right] = 0
\end{equation*}
for any $n \geq 1$, any $\{a_1,\dots,a_n,b_1,\dots,b_n\} \subset \reels$ and any (sufficiently integrable) functions $\{h_1,\dots,h_n,g_1,\dots,g_n\}$. Then, for a fixed $n$ and fixed functions $\{h_1,\dots,h_n,g_1,\dots,g_n\}$, the Laplace transform of the signed measure on $\mathcal{B}(\reels^n) \times \mathcal{B}(\reels^n)$ defined by
\begin{equation*}
(A,B) \mapsto \esp \left[ X \ind_A \bigr( W(h_1), \dots, W(h_n) \bigr) \ind_B \bigr( \widetilde{N}(g_1), \dots, \widetilde{N}(g_n) \bigr) \right] ,
\end{equation*}
is identically $0$. Consequently, the measure on $\tribu = \tribu_T$ defined by $E \mapsto \esp \left[ X \ind_E \right]$ vanishes on every rectangle $A \times B$ if it is a pre-image of the $\reels^{2n}$-dimensional random vector
$$
\left(W(h_1), \dots, W(h_n), \widetilde{N}(g_1), \dots, \widetilde{N}(g_n) \right) .
$$
By linearity of the stochastic integrals, this is also true for random vectors of the form
$$
\left(W(h_1), \dots, W(h_n), \widetilde{N}(g_1), \dots, \widetilde{N}(g_m) \right) ,
$$
when $m$ and $n$ are different. Since $\tribu$ is generated by those random vectors, the measure is identically zero and $X = 0$.
\end{proof}

We now state and prove a Martingale Representation Theorem with respect to the Brownian motion and the Poisson random measure simultaneously.
\begin{thm}\label{T:prp}
Let $F \in L^2(\Omega, \tribu, \p)$. There exist a unique Borel predictable process $\psi \in L^2(\mesleb \times \meslevy \times \p)$ and a unique predictable process $\phi \in L^2(\mesleb \times \p)$ such that
\begin{equation}\label{E:repLevy}
F = \esp[F] + \int_0^T \phi(t) \, W(dt) + \intPoissonT \psi(t,z) \, \Ntildedt .
\end{equation}
\end{thm}
\begin{proof}
For $h \in L^2([0,T], \mesleb)$ and $g \in C([0,T])$, we know from the proof of Lemma~\ref{L:procZ} that
\begin{multline*}
Y_t = \exp \left\lbrace \int_0^t h(s) \, W(ds) - \frac{1}{2} \int_0^t h^2(s) \, ds + \intPoissont g(s) \gamma(z) \, \Ntildeds \right.\\
\left. - \intPoissont \left( e^{g(s) \gamma(z)} - 1 - g(s) \gamma(z) \right) \, \meslevy(dz) ds \right\rbrace
\end{multline*}
is a solution of
\begin{equation}
Y_t = 1 + \int_0^t Y_{s-} h(s) \, W(ds) + \int_0^t \int_{\reels} Y_{s-} \left( e^{g(s) \gamma(z)} - 1 \right) \, \Ntildeds
\end{equation}
on $[0,T]$. Hence, $Y_T$ admits a martingale representation as in Equation~\eqref{E:repLevy} with $\phi(t) = Y_{t-} h(t)$ and $\psi(t,z) = Y_{t-} \left( e^{g(t) \gamma(z)} - 1 \right)$. These two processes are predictable. Note that
$$
Y_T = Y(h,g) e^{- \theta_T (h,g)}
$$
where
$$
\theta_T (h,g) = \frac{1}{2} \int_0^T h^2(t) \, dt + \intPoissonT \left( e^{g(t) \gamma(z)} - 1 - g(t) \gamma(z) \right) \, \meslevy(dz) dt .
$$
Since $\theta_T (h,g)$ is deterministic, $Y(h,g)$ also admits a martingale representation as in Equation~\eqref{E:repLevy} but this time with
$$
\phi(t) = Y_{t-} h(t) e^{ \theta_T (h,g)} \quad \text{and} \quad \psi(t,z) = Y_{t-} \left( e^{g(t) \gamma(z)} - 1 \right) e^{\theta_T (h,g)} .
$$
Therefore, the first statement follows by a denseness argument. Indeed, from Lemma~\ref{L:densite}, since $F$ is square-integrable, there exists a sequence $(F_n)_{n \geq 1}$ of square-integrable random variables such that $F_n$ tends to $F$ in the $L^2(\Omega)$-norm when $n$ goes to infinity. Moreover, the $F_n$'s are linear combinations of some $Y(h,g)$'s. Then, for each term in this sequence there exist $\phi_n$ and $\psi_n$ such that
$$
F_n = \esp[F_n] + \int_0^T \phi_n (t) \, W(dt) + \intPoissonT \psi_n (t,z) \, \Ntildedt .
$$
Also, since
\begin{align*}
\esp [F_n - F_m]^2 &= \esp \biggr[ \esp[F_n - F_m] + \int_0^T (\phi_n (t) - \phi_m (t)) \, W(dt) \\
& \qquad \qquad + \intPoissonT (\psi_n (t,z) - \psi_m (t,z)) \, \Ntildedt \biggr]^2 \\
&= \left( \esp[F_n - F_m]\right)^2 + \int_0^T \esp [\phi_n (t) - \phi_m (t)]^2 \, dt \\
& \qquad \qquad + \intPoissonT \esp [\psi_n (t,z) - \psi_m (t,z)]^2 \, \meslevy(dz) dt ,
\end{align*}
we get that $(\phi_n)_{n \geq 1}$ and $(\psi_n)_{n \geq 1}$ are Cauchy sequences. It follows that there exist predictable processes $\psi \in L^2(\mesleb \times \meslevy \times \p)$ and $\phi \in L^2(\mesleb \times \p)$ for which the representation of Equation~\eqref{E:repLevy} is verified.

We now prove the second statement. If $F$ admits two martingale representations with $\phi_1, \phi_2, \psi_1, \psi_2$ (these have nothing to do with the previous sequences), then by Itô's isometry
$$0 = \| \phi_1 - \phi_2 \|^2_{L^2(\mesleb \times \p)} + \| \psi_1 - \psi_2 \|^2_{L^2(\mesleb \times \meslevy \times \p)}$$
and then $\phi_1 = \phi_2$ in $L^2(\mesleb \times \p)$ and $\psi_1 = \psi_2$ in $L^2(\mesleb \times \meslevy \times \p)$.
\end{proof}

\begin{rem}
From now on, we will refer to this martingale representation property of the Lévy process $X$ as the MRP.
\end{rem}

\section{Chaotic representations}

We now define multiple integrals with respect to $W(dt)$ and $\Ntildedt$ simultaneously and define Lévy chaos as an extension of Wiener-Itô chaos. Then, we show that any square-integrable Lévy functional can be represented by a chaos expansion. We refer the reader to the lecture notes of Meyer \cite{meyer1976} for more details on multiple stochastic integrals.

\subsection{Notation}

In the following, we unify the notation of the Poisson random measure and the Brownian motion. Thus, the superscript $(1)$ will refer to Brownian motion and the superscript $(2)$ to the Poisson random measure. This is also the notation in \cite{benthetal2003}.

Let $\x = [0,T] \times \reels$. We introduce two (projection) operators $\Pi_1 \colon \x \to [0,T]$ and $\Pi_2 \colon \x \to \x$ defined by $\Pi_1(t,z) = t$ and $\Pi_2(t,z) = (t,z)$. Consequently, $\Pi_1 \left( [0,T] \times \reels \right) = [0,T]$ and $\Pi_2 \left( [0,T] \times \reels \right) = [0,T] \times \reels$.

For $n \geq 1$, $t \in [0,T]$ and $(i_1,\dots,i_n) \in \{1,2\}^n$, we also introduce the following notations:
\begin{equation}\label{E:defsigmant}
\Sigma_n (t) = \left\lbrace (t_1,\dots,t_n) \in [0,T]^n \mid t_1 < \dots < t_n \leq t \right\rbrace ;
\end{equation}
and
\begin{multline*}
\Sigma_{(i_1,\dots,i_n)} ([0,t] \times \reels) \\
= \left\lbrace (x_1,\dots,x_n) \in \Pi_{i_1}(\x) \times \dots \times \Pi_{i_n}(\x) \mid \Pi_1 (x_1) < \dots < \Pi_1 (x_n) \leq t \right\rbrace .
\end{multline*}
Consequently, $\Sigma_n (T) = \Sigma_{(i_1,\dots,i_n)} (\x)$ when $i_k = 1$ for each $k = 1,2,\dots,n$. If $f$ is a function defined on $\zerotn$, we write $f(x_1,\dots,x_n)$, where $x_k \in \Pi_{i_k}(\x)$ for each $k = 1,2,\dots,n$. If $\eta_1 = \mesleb$ and $\eta_2 = \mesleb \times \meslevy$, let $L^2 \left( \sigman \right)$ be the space of square-integrable functions defined on $\sigman$ and equipped with the product measure $\eta_{i_1} \times \dots \times \eta_{i_n}$ defined on $\Pi_{i_1}(\x) \times \dots \times \Pi_{i_n}(\x)$.

\subsection{Multiple integrals and Lévy chaos}

Fix $n \geq 1$ and $(i_1,\dots,i_n) \in \{1,2\}^n$. We define the iterated integral $J_{(i_1,\dots,i_n)} (f)$, for $f$ in $L^2 \left( \sigman \right)$, by
\begin{multline*}
J_{(i_1,\dots,i_n)} (f) \\
= \int_{\Pi_{i_n}([0,T] \times \reels)} \dots \int_{\Pi_{i_1}([0,t_2-] \times \reels)} f(x_1,\dots,x_n) \, M^{(i_1)}(dx_1) \dots M^{(i_n)}(dx_n)
\end{multline*}
where $M^{(j)}(dx)$ equals $W(dt)$ if $j = 1$ and equals $\Ntildedt$ if $j = 2$. The $i_1$ in $J_{(i_1,\dots,i_n)}$ stands for the innermost stochastic integral and the $i_n$ stands for the outermost stochastic integral. For example, if $n = 3$ and $(i_1,i_2,i_3) = (1,1,2)$, then
\begin{multline*}
J_{(1,1,2)} (f) \\
= \int_0^T \int_{\reels} \left[ \int_0^{t_3-} \left( \int_0^{t_2-} f(t_1,t_2,(t_3,z_3)) W(dt_1) \right) W(dt_2) \right] \widetilde{N}(dt_3,dz_3) .
\end{multline*}

As $n$ runs through $\nat$ and $(i_1,\dots,i_n)$ runs through $\{1,2\}^n$, the iterated integrals generate orthogonal spaces in $L^2(\Omega)$ that we would like to call \textit{Lévy chaos}. Indeed, since
$$
\int_{\Pi_{i}([0,T] \times \reels)} f (x) M^{(i)}(dx)
$$
and
$$
\int_{\Pi_{j}([0,T] \times \reels)} g (x) M^{(j)}(dx)
$$
are independent if $i \neq j$ and both have mean zero, using Itô's isometry iteratively, we get the following proposition.

\begin{prop}\label{P:orthogonality}
If $f \in L^2(\sigman)$ and if $g \in L^2(\sigmam)$, then
\begin{multline*}
\esp \left[ J_{(i_1,\dots,i_n)}(f) J_{(j_1,\dots,j_m)}(g)\right] \\
=
\begin{cases}
(f,g)_{L^2(\sigman)} & \text{if $(i_1,\dots,i_n) = (j_1,\dots,j_m)$;}\\
0 & \text{if not.}
\end{cases}
\end{multline*}
\end{prop}

We end this subsection with a definition.

\begin{defn}\label{D:tensorproduct}
For $n \geq 1$ and $(i_1,\dots,i_n) \in \{1,2\}^n$, the $(i_1,\dots,i_n)$-tensor product of a function $h$ defined on $[0,T]$ with a function $g$ defined on $[0,T] \times \reels$ is a function on $\zerotn$ defined by
$$
\left( h \otimes_{(i_1,\dots,i_n)} g \right) (x_1,\dots,x_n) = \prod_{1 \leq k \leq n} h \left( \Pi_1(x_k) \right)^{2 - i_k} g \left( \Pi_2(x_k) \right)^{i_k - 1} .
$$
\end{defn}

For example,
$$
\left( h \otimes_{(1,1)} g \right) (s,t) = h(s)h(t)
$$
is a function defined on $[0,T] \times [0,T]$ and
$$
\left( h \otimes_{(1,2,1)} g \right) (r,(s,y),t) = h(r)h(t)g(s,y)
$$
is a function defined on $[0,T] \times ([0,T] \times \reels) \times [0,T]$.

\subsection{Chaotic representation property}

For the rest of the paper, we will assume that $\sum_{(i_1,\dots,i_n)}$ means $\sum_{(i_1,\dots,i_n) \in \{1,2\}^n}$.

Recall that $Z = (Z_t)_{t \in [0,T]}$ was defined in Equation~\eqref{E:procZ} by
\begin{multline*}
Z_t = \exp \left\lbrace \int_0^t h(s) \, W(ds) - \frac{1}{2} \int_0^t h^2(s) \, ds + \intPoissont g(s,z) \, \Nds \right. \\
- \left. \intPoissont \left( e^{g(s,z)} - 1 \right) \, \meslevy(dz) ds \right\rbrace .
\end{multline*}

\begin{lem}\label{L:crp}
Let $h \in L^2([0,T])$ and $e^g - 1 \in L^2([0,T] \times \reels, \mesleb \times \meslevy)$. Then, $Z_T$ admits the following chaotic representation:
\begin{equation}
Z_T = 1 + \sum_{n = 1}^{\infty} \sum_{(i_1,\dots,i_n)} J_{(i_1,\dots,i_n)} \left( h \otimes_{(i_1,\dots,i_n)} (e^{g} - 1) \right) .
\end{equation}
\end{lem}
\begin{proof}
We know from the proof of Lemma~\ref{L:procZ} that $Z_T$ is square-integrable and that
\begin{equation}\label{E:temp}
Z_T = 1 + \int_0^T Z_{t-} h(t) \, W(dt) + \intPoissonT Z_{t-} (e^{g(t,z)} - 1) \, \Ntildedt .
\end{equation}
Let $\phi^{(1)} (t) = Z_{t-} h(t)$ and $\phi^{(2)} (t,z) = Z_{t-} (e^{g(t,z)} - 1)$. We now iterate Equation~\eqref{E:temp}. Consequently,
\begin{align*}
Z_T = 1 &+ \int_0^T f^{(1)}(t) \, W(dt) + \intPoissonT f^{(2)}(t,z) \, \Ntildedt \\
&+ \int_0^T \int_0^{t-} Z_{s-} h(s) h(t) \, W(ds) \, W(dt) \\
&+ \int_0^T \intPoissontmoins Z_{s-} (e^{g(s,y)} - 1) h(t) \, \Ntildedyds \, W(dt) \\
&+ \intPoissonT \int_0^{t-} Z_{s-} h(s) (e^{g(t,z)} - 1) \, W(ds) \, \Ntildedt \\
&+ \intPoissonT \intPoissontmoins Z_{s-} (e^{g(s,y)} - 1) (e^{g(t,z)} - 1) \, \Ntildedyds \, \Ntildedt
\end{align*}
where $f^{(1)}(t) = h(t) = \left( h \otimes_{(1)} (e^g - 1) \right) (t)$ and $f^{(2)}(t,z) = e^{g(t,z)} - 1 = \left( h \otimes_{(2)} (e^g - 1) \right) (t,z)$. Then, after $n$ iterations, we get
\begin{multline*}
Z_T = 1 + \sum_{k = 1}^{n-1} \sum_{(i_1,\dots,i_k)} J_{(i_1,\dots,i_k)} (f^{(i_1,\dots,i_k)}) \\
+ \sum_{(i_1,\dots,i_n)} \int_{\Pi_{i_n}([0,T] \times \reels)} \dots \int_{\Pi_{i_1}([0,t_2-] \times \reels)} \phi^{(i_1,\dots,i_n)} \left( x_1,\dots,x_n \right) \\
M^{(i_1)}(dx_1) \dots M^{(i_n)}(dx_n)
\end{multline*}
where $f^{(i_1,\dots,i_k)} = h \otimes_{(i_1,\dots,i_k)} (e^{g} - 1)$ and where $\phi^{(i_1,\dots,i_n)} = Z_{-} (h \otimes_{(i_1,\dots,i_n)} (e^{g} - 1))$. This means that we can define a sequence $(\psi_n)_{n \geq 2}$ in $L^2(\Omega)$ by
\begin{multline*}
\psi_n = \sum_{(i_1,\dots,i_n)} \int_{\Pi_{i_n}([0,T] \times \reels)} \dots \int_{\Pi_{i_1}([0,t_2-] \times \reels)} \phi^{(i_1,\dots,i_n)} \left( x_1,\dots,x_n \right) \\
M^{(i_1)}(dx_1) \dots M^{(i_n)}(dx_n) .
\end{multline*}

From Proposition~\ref{P:orthogonality},
$$
\esp[Z_T^2] = 1 + \sum_{k = 1}^{n-1} \sum_{(i_1,\dots,i_k)} \| f^{(i_1,\dots,i_k)} \|_{L^2(\Sigma_{(i_1,\dots,i_k)}(\x))}^2 + \esp[\psi_n^2]
$$
for each $n \geq 2$. Hence we get that
$$
\sum_{n = 1}^{\infty} \sum_{(i_1,\dots,i_n)} J_{(i_1,\dots,i_n)} (f^{(i_1,\dots,i_n)})
$$
is a square-integrable series and that there exists a square-integrable random variable $\psi$ such that $\psi_n$ tends to $\psi$ in the $L^2(\Omega)$-norm. Consequently, it is enough to show that $\psi = 0$. Since $f^{(i_1,\dots,i_n)} = h \otimes_{(i_1,\dots,i_n)} (e^{g} - 1)$, using Proposition~\ref{P:orthogonality} once again, we get that
\begin{multline*}
\sum_{(i_1,\dots,i_n)} \esp \left[ \left( J_{(i_1,\dots,i_n)} (f^{(i_1,\dots,i_n)}) \right)^2 \right] \\
= \sum_{k = 0}^n \sum_{\stackrel{(i_1,\dots,i_n)}{|i|=k}} \| h \otimes_{(i_1,\dots,i_n)} (e^{g} - 1) \|^2_{L^2(\sigman)} ,
\end{multline*}
where $|i| = |(i_1,\dots,i_n)| = \sum_{j = 1}^n (2-i_j)$ stands for the number of times the function $h$ appears in the tensor product. Note that when $|i| = k$ there are $\binom{n}{k}$ terms in the innermost summation. Since $h^2 \otimes_{(i_1,\dots,i_n)} (e^{g} - 1)^2$ is a
$(i_1,\dots,i_n)$-tensor product, the function given by
$$
\sum_{\stackrel{(i_1,\dots,i_n)}{|i|=k}} h \otimes_{(i_1,\dots,i_n)} (e^{g} - 1)
$$
is symmetric on $\zerotn$. Consequently,
\begin{align*}
& \sum_{(i_1,\dots,i_n)} \esp \left[ \left( J_{(i_1,\dots,i_n)} (f^{(i_1,\dots,i_n)}) \right)^2 \right] \\
&= \sum_{k = 0}^n \int_{\sigman} \biggr[ \sum_{\stackrel{(i_1,\dots,i_n)}{|i|=k}} h^2 \otimes_{(i_1,\dots,i_n)} (e^{g} - 1)^2 \biggr] \, d\eta_{i_1} \dots d\eta_{i_n}\\
&= \frac{1}{n!} \sum_{k = 0}^n \int_{\zerotn} \biggr[ \sum_{\stackrel{(i_1,\dots,i_n)}{|i|=k}} h^2 \otimes_{(i_1,\dots,i_n)} (e^{g} - 1)^2 \biggr] \, d\eta_{i_1} \dots d\eta_{i_n}\\
&= \frac{1}{n!} \sum_{k = 0}^n \binom{n}{k} \| h \|^{2k}_{L^2(\mesleb)} \|e^{g} - 1\|^{2(n-k)}_{L^2(\mesleb \times \meslevy)} \\
&= \frac{1}{n!} \left( \| h \|^2_{L^2(\mesleb)} + \|e^{g} - 1\|^2_{L^2(\mesleb \times \meslevy)} \right)^n .
\end{align*}

From Equation~\eqref{E:lemprocZ}, we know that
$$
\esp [Z_T^2] = \exp \left\lbrace \| h \|^2_{L^2(\mesleb)} + \|e^{g} - 1\|^2_{L^2(\mesleb \times \meslevy)} \right\rbrace .
$$
This means that $\psi = 0$ and the statement follows.
\end{proof}

We are now ready to state and prove the chaotic representation property of the Lévy process $X$. The previous lemma and the idea of its proof will be of great use.

\begin{thm}\label{T:crp}
Let $F \in L^2(\Omega, \tribu, \p)$. There exists a unique sequence
$$
\left\lbrace f^{(i_1,\dots,i_n)}; n \geq 1, (i_1,\dots,i_n) \in \{1,2\}^n \right\rbrace ,
$$
whose elements are respectively in $L^2 \left( \sigman \right)$, such that
\begin{equation}\label{E:chaos}
F = \esp[F] + \sum_{n = 1}^{\infty} \sum_{(i_1,\dots,i_n)} J_{(i_1,\dots,i_n)} \left( f^{(i_1,\dots,i_n)} \right) .
\end{equation}
Consequently,
\begin{equation}\label{E:isometrie}
\esp[F^2] = \esp^2[F] + \sum_{n = 1}^{\infty} \sum_{(i_1,\dots,i_n)} \| f^{(i_1,\dots,i_n)} \|_{L^2(\sigman)}^2 .
\end{equation}
\end{thm}
\begin{proof}
From Theorem~\ref{T:prp}, we know there exist a predictable process $\phi^{(1)} \in L^2(\mesleb \times \p)$ and a Borel predictable process $\phi^{(2)} \in L^2(\mesleb \times \meslevy \times \p)$ such that
\begin{equation*}
F = \esp[F] + \int_0^T \phi^{(1)} (t) \, W(dt) + \intPoissonT \phi^{(2)} (t,z) \, \Ntildedt .
\end{equation*}
Using Itô's isometry, it is clear that
$$
\| \phi^{(1)} \|_{L^2(\mesleb \times \p)}^2 + \| \phi^{(2)} \|_{L^2(\mesleb \times \meslevy \times \p)}^2 \leq \esp[F^2] .
$$
For almost all $t \in [0,T]$, $\phi^{(1)}(t) \in L^2(\Omega, \tribu_t, \p)$ and then from Theorem~\ref{T:prp} there exist processes $\phi^{(1,1)}$ and $\phi^{(1,2)}$ such that
\begin{equation*}
\phi^{(1)}(t) = \esp[\phi^{(1)}(t)] + \int_0^t \phi^{(1,1)} (t,s) \, W(ds) + \intPoissont \phi^{(1,2)} (t,s,y) \, \Ntildedyds .
\end{equation*}
Similarly, for almost all $(t,z) \in [0,T] \times \reels$, $\phi^{(2)}(t,z) \in L^2(\Omega, \tribu_t, \p)$ and
\begin{multline*}
\phi^{(2)}(t,z) = \esp[\phi^{(2)}(t,z)] + \int_0^t \phi^{(2,1)} (t,z,s) \, W(ds) \\
+ \intPoissont \phi^{(2,2)} (t,z,s,y) \, \Ntildedyds .
\end{multline*}
Consequently,
\begin{multline*}
F = \esp[F] + \int_0^T g^{(1)}(t) \, W(dt) + \intPoissonT g^{(2)}(t,z) \, \Ntildedt \\
+ \int_0^T \int_0^{t-} \phi^{(1,1)} (t,s) \, W(ds) \, W(dt) \\
+ \int_0^T \intPoissontmoins \phi^{(1,2)} (t,s,y) \, \Ntildedyds \, W(dt) \\
+ \intPoissonT \int_0^{t-} \phi^{(2,1)} (t,z,s) \, W(ds) \, \Ntildedt \\
+ \intPoissonT \intPoissontmoins \phi^{(2,2)} (t,z,s,y) \, \Ntildedyds \, \Ntildedt .
\end{multline*}
where $g^{(1)}(t) = \esp[\phi^{(1)}(t)]$ and $g^{(2)}(t,z) = \esp[\phi^{(2)}(t,z)]$. After $n$ steps of this procedure, i.e. after $n$ iterations of Theorem~\ref{T:prp}, we get as in the proof of Lemma~\ref{L:crp} that
$$
F = \esp[F] + \sum_{k = 1}^{n-1} \sum_{(i_1,\dots,i_k)} J_{(i_1,\dots,i_k)} (f^{(i_1,\dots,i_k)}) + \psi_n
$$
where $f^{(i_1,\dots,i_k)} \in L^2 \left( \Sigma_{(i_1,\dots,i_k)}(\x) \right)$, for each $1 \leq k \leq n-1$ and $(i_1,\dots,i_k) \in \{1,2\}^k$, where
\begin{multline*}
\psi_n = \sum_{(i_1,\dots,i_n)} \int_{\Pi_{i_n}([0,T] \times \reels)} \dots \int_{\Pi_{i_1}([0,t_2-] \times \reels)} \phi^{(i_1,\dots,i_n)} \left( x_1,\dots,x_n \right) \\
M^{(i_1)}(dx_1) \dots M^{(i_n)}(dx_n) ,
\end{multline*}
and where $\phi^{(i_1,\dots,i_n)} \in L^2 \left( \eta_{i_1} \times \dots \times \eta_{i_n} \times \p \right)$, for each $(i_1,\dots,i_n) \in \{1,2\}^n$.

From Proposition~\ref{P:orthogonality},
$$
\esp[F^2] = \esp[F]^2 + \sum_{k = 1}^{n-1} \sum_{(i_1,\dots,i_k)} \| f^{(i_1,\dots,i_k)} \|_{L^2(\Sigma_{(i_1,\dots,i_k)}(\x))}^2 + \esp[\psi_n^2] ,
$$
for each $n \geq 2$ and
$$
\sum_{n = 1}^{\infty} \sum_{(i_1,\dots,i_n)} J_{(i_1,\dots,i_n)} (f^{(i_1,\dots,i_n)})
$$
is a square-integrable series. Consequently, we know that there exists a square-integrable random variable $\psi$ such that $\psi_n$ tends to $\psi$ in the $L^2(\Omega)$-norm. It is enough to show that $\psi = 0$. Using the argument leading to Proposition~\ref{P:orthogonality}, i.e. the fact that two iterated stochastic integrals of different order are orthogonal, we get that for a fixed $n \geq 2$,
$$
\left( J_{(i_1,\dots,i_k)} (f^{(i_1,\dots,i_k)}), \psi_n \right)_{L^2(\Omega)} = 0
$$
for every $1 \leq k \leq n-1$, $(i_1,\dots,i_k) \in \{1,2\}^k$ and $f^{(i_1,\dots,i_k)} \in L^2(\Sigma_{(i_1,\dots,i_k)}(\x))$. Thus,
\begin{equation}\label{E:psiortho}
\left( J_{(i_1,\dots,i_n)} (f^{(i_1,\dots,i_n)}), \psi \right)_{L^2(\Omega)} = 0
\end{equation}
for every $n \geq 1$, $(i_1,\dots,i_n) \in \{1,2\}^n$ and $f^{(i_1,\dots,i_n)} \in L^2(\sigman)$.

We now assume that $g = \bar{g} \gamma$ where $\bar{g}$ belongs to $C([0,T])$. Using Equation~\eqref{E:psiortho}, we have that $\psi$ is orthogonal to each random variable $Y (h,g)$ defined in Equation~\eqref{E:dense} since from Lemma~\ref{T:crp} they each possess a chaos decomposition. We also know from Lemma~\ref{L:densite} that these random variables are dense in $L^2(\Omega, \tribu, \p)$, so $\psi = 0$. This means that every square-integrable Lévy functional can be express as a series of iterated integrals. The statement follows.
\end{proof}

\begin{rem}
From now on, we will refer to the chaotic representation property of Theorem~\ref{T:crp} as the CRP.
\end{rem}

\begin{rem}\label{R:crpimpliesmrp}
As mentioned before, in general the CRP implies the MRP. Indeed, if $F$ is a square-integrable Lévy functional with chaos decomposition
$$
F = \esp[F] + \sum_{n = 1}^{\infty} \sum_{(i_1,\dots,i_n)} J_{(i_1,\dots,i_n)} \left( f^{(i_1,\dots,i_n)} \right) ,
$$
then
$$
F = \esp[F] + \int_0^T \phi(t) \, W(dt) + \intPoissonT \psi(t,z) \, \Ntildedt ,
$$
with
\begin{equation*}
\begin{aligned}
\phi(t) &= f^{(1)}(t) + \sum_{n = 1}^{\infty} \sum_{(i_1,\dots,i_n)} J_{(i_1,\dots,i_n)} \left( f^{(i_1,\dots,i_n,1)} (\cdot,t) \ind_{\Sigma_n (t)} \right) ,\\
\psi(t,z) &= f^{(2)}(t,z) + \sum_{n = 1}^{\infty} \sum_{(i_1,\dots,i_n)} J_{(i_1,\dots,i_n)} \left( f^{(i_1,\dots,i_n,2)} (\cdot,(t,z)) \ind_{\Sigma_n (t)} \right) .
\end{aligned}
\end{equation*}
\end{rem}

This last remark, together with our journey from the MRP of Theorem~\ref{T:prp} to the CRP of Theorem~\ref{T:crp}, yields the following interesting proposition.

\begin{prop}
For a square-integrable Lévy process, the MRP and the CRP are equivalent.
\end{prop}

\subsection{Explicit chaos representation}

In the next proposition, we compute the explicit chaos representation of a \textit{smooth} Lévy functional.

\begin{prop}
Let $f$ be a smooth function with compact support in $\reels^k$, i.e. let $f \in C^{\infty}_c (\reels^k)$, and let $t_j$ belong to $[0,T]$ for each $j=1,\dots,k$. Then,
$$
f(X_{t_1},\dots,X_{t_k}) = \esp[f(X_{t_1},\dots,X_{t_k})] + \sum_{n = 1}^{\infty} \sum_{(i_1,\dots,i_n)} J_{(i_1,\dots,i_n)} (f^{(i_1,\dots,i_n)}) ,
$$
where
\begin{multline*}
f^{(i_1,\dots,i_n)}(\Pi_{i_1}(s_1,w_1),\dots,\Pi_{i_n}(s_n,w_n)) \\
=  - (2 \pi)^{-k/2} \int_{\reels^k} \hat{f}(y) \widehat{\phi} (-y) \prod_{1 \leq j \leq n} (i \sigma \xi^{t,y}_{s_j})^{2 - i_j} (e^{i w_j \xi^{t,y}_{s_j}} - 1)^{i_j - 1} \, dy .
\end{multline*}
with
$$
\phi(x) \, dx = \p\{X_t \in dx\} ,
$$
where $X_t = (X_{t_1},\dots,X_{t_k})$, and with
$$
\xi_s^{t,y} = y_1 \ind_{[0,t_1]}(s) + \dots + y_k \ind_{[0,t_k]}(s) ,
$$
for $t = (t_1,\dots,t_k)$ and $y = (y_1,\dots,y_k)$.
\end{prop}
\begin{proof}
We follow an idea in \cite{lokka2004} and use Fourier transforms. By the Fourier inversion formula
\begin{equation}\label{E:fourierinvf}
f(x) = (2 \pi)^{-k/2} \int_{\reels^k} \hat{f}(y) e^{i \langle x,y \rangle} \, dy
\end{equation}
where $\hat{f}$ is the Fourier transform of $f$ and $\langle x,y \rangle$ denotes the scalar product in $\reels^k$ of $x = (x_1,\dots,x_k)$ and $y = (y_1,\dots,y_k)$. Let $\phi^- (x) = \phi(-x)$. If we define $F(x) = \esp \left[ f(X_{t_1} + x_1,\dots,X_{t_k} + x_k) \right]$, then
$$
F(x) = -(f \ast \phi^-)(x)
$$
and also
$$
F(x) = -(2 \pi)^{-k/2} \int_{\reels^k} \hat{f}(y) \widehat{\phi^-}(y) e^{i \langle x,y \rangle} \, dy .
$$
Therefore, we have the following equality:
$$
\esp[f(X_{t_1},\dots,X_{t_k})] = -(2 \pi)^{-k/2} \int_{\reels^k} \hat{f}(y) \widehat{\phi^-}(y) \, dy .
$$

From Equation~\eqref{E:fourierinvf} and Equation~\eqref{E:Lgen}, we have that
\begin{align*}
f(X_{t_1},\dots,X_{t_k}) &= (2 \pi)^{-k/2} \int_{\reels^k} \hat{f}(y) e^{i \langle X_t,y \rangle } \, dy \\
&= (2 \pi)^{-k/2} \int_{\reels^k} \hat{f}(y) e^{i \mu \langle t,y \rangle } Y^{t,y} \, dy
\end{align*}
where
$$
Y^{t,y} = \exp \left\lbrace \int_0^T i \sigma \xi_s^{t,y} \, W(ds) + \intPoissonT i z \xi_s^{t,y} \, \Ntildeds \right\rbrace .
$$
Hence,
\begin{equation}\label{E:fdeX}
f(X_{t_1},\dots,X_{t_k}) = (2 \pi)^{-k/2} \int_{\reels^k} \hat{f}(y) e^{i \mu \langle t,y \rangle } Z^{t,y} \esp [Y^{t,y}] \, dy
\end{equation}
where
\begin{multline*}
Z^{t,y} = \exp \biggr\lbrace \int_0^T i \sigma \xi_s^{t,y} \, W(ds) + \frac{1}{2} \sigma^2 \int_0^T (\xi_s^{t,y})^2 \, ds \\
+ \intPoissonT i z \xi_s^{t,y} \, \Nds - \intPoissonT (e^{i z \xi_s^{t,y}} - 1) \, \meslevy(dz) ds \biggr\rbrace .
\end{multline*}
From Lemma~\ref{L:crp}, we know that
$$
Z^{t,y} = 1 +  \sum_{n = 1}^{\infty} \sum_{(i_1,\dots,i_n)} J_{(i_1,\dots,i_n)} \left( (i \sigma \xi^{t,y}) \otimes_{(i_1,\dots,i_n)} (e^{i z \xi^{t,y}} - 1) \right) .
$$

On the other hand,
\begin{align*}
\esp [Y^{t,y}] &= e^{-i \mu \langle t,y \rangle } \esp \left[  e^{i \langle X_t,y \rangle } \right] \\
&= - e^{-i \mu \langle t,y \rangle } \widehat{\phi} (-y)
\end{align*}
Then, using Equation~\eqref{E:fdeX} and by Lebesgue's dominated convergence theorem,
\begin{align*}
f(& X_{t_1},\dots,X_{t_k}) \\
&= - (2 \pi)^{-k/2} \int_{\reels^k} \hat{f}(y) \widehat{\phi} (-y) \, dy \\
& \quad \quad - (2 \pi)^{-k/2} \int_{\reels^k} \hat{f}(y) \widehat{\phi} (-y) \\
& \qquad \qquad \qquad \times \sum_{n = 1}^{\infty} \sum_{(i_1,\dots,i_n)} J_{(i_1,\dots,i_n)} \left( (i \sigma \xi^{t,y}) \otimes_{(i_1,\dots,i_n)} (e^{i z \xi^{t,y}} - 1) \right) \, dy \\
&= \esp[f(X_{t_1},\dots,X_{t_k})] + \sum_{n = 1}^{\infty} \sum_{(i_1,\dots,i_n)} J_{(i_1,\dots,i_n)} (f^{(i_1,\dots,i_n)}) .
\end{align*}
where
\begin{multline*}
f^{(i_1,\dots,i_n)}(\Pi_{i_1}(s_1,w_1),\dots,\Pi_{i_n}(s_n,w_n)) \\
=  - (2 \pi)^{-k/2} \int_{\reels^k} \hat{f}(y) \widehat{\phi} (-y) \left( (i \sigma \xi^{t,y}) \otimes_{(i_1,\dots,i_n)} (e^{i z \xi^{t,y}} - 1) \right) \, dy .
\end{multline*}
The statement follows from Definition~\ref{D:tensorproduct}.
\end{proof}

\section{Malliavin derivatives and Clark-Ocone formula}

Before defining the Malliavin derivatives, we introduce a last notation: for $n \geq 1$ and $1 \leq k \leq n+1$, define
$$
\Sigma^k_n (t) = \left\lbrace (t_1,\dots,t_n) \in [0,T]^n \mid t_1 < \dots < t_{k-1} < t < t_k < \dots < t_n \right\rbrace ,
$$
i.e. $t$ is at the $k$-th position between the $t_j$'s, where $t_0 = 0$ and $t_{n+1} = T$. Note that $\Sigma^{n+1}_n (t) = \Sigma_n (t)$, where the latter was defined earlier in Equation~\eqref{E:defsigmant}. In a multi-index $(i_1,\dots,i_n)$, we will use $\widehat{i}_k$ to denote the omission of the $k$-th index.

We want to define two directional derivative operators in the spirit of Le{\'o}n et al. \cite{leonetal2002}: one in the direction of the Brownian motion and one in the direction of the Poisson random measure. If $F = J_{(i_1,\dots,i_n)} (f)$, then we would like to define $D^{(1)}_t F$ and $D^{(2)}_{t,z} F$ as follows:
$$
D^{(1)}_t F = \sum_{k=1}^n \ind_{\{i_k = 1\}}  J_{(i_1,\dots,\widehat{i}_k,\dots,i_n)} \bigr( f(\underbrace{\dots}_{k-1},t,\underbrace{\dots}_{n - k}) \ind_{\Sigma^k_{n-1} (t)} \bigr)
$$
and
$$
D^{(2)}_{t,z} F = \sum_{k=1}^n \ind_{\{i_k = 2\}}  J_{(i_1,\dots,\widehat{i}_k,\dots,i_n)} \bigr( f(\underbrace{\dots}_{k-1},(t,z),\underbrace{\dots}_{n - k}) \ind_{\Sigma^k_{n-1} (t)} \bigr)
$$
where $J_{(\, \widehat{i} \,)} (f) = f$.

\begin{defn}
Let $\D = \mathbb{D}^{(1)} \cap \mathbb{D}^{(2)}$, where if $j = 1$ or if $j = 2$, $\mathbb{D}^{(j)}$ is the subset of $L^2(\Omega, \tribu, \p)$ consisting of the random variables $F$ with chaotic representation
$$
F = \esp[F] + \sum_{n = 1}^{\infty} \sum_{(i_1,\dots,i_n)} J_{(i_1,\dots,i_n)} \bigr( f^{(i_1,\dots,i_n)} \bigr)
$$
such that
$$
\sum_{n = 1}^{\infty} \sum_{(i_1,\dots,i_n)} \sum_{k=1}^n \ind_{\{i_k = j\}} \int_{\Pi_{j}(\x)} \bigr\| f^{(i_1,\dots,i_n)} (\cdot,x,\cdot) \ind_{\Sigma^k_{n-1} (t)} \bigr\|^2 \, \eta_{j}(dx) < \infty ,
$$
where the inside norm is the $L^2(\sigmanomisk)$-norm.
\end{defn}

From Theorem~\ref{T:crp}, it is clear that $\D$ is dense in $L^2(\Omega)$, since every random variable with a chaos representation given by a finite sum belongs to $\D$.

\begin{defn}\label{D:derivees}
The Malliavin derivatives $D^{(1)} \colon \mathbb{D}^{(1)} \to L^2 \left( [0,T] \times \Omega \right)$ and $D^{(2)} \colon \mathbb{D}^{(2)} \to L^2 \left( [0,T] \times \reels \times \Omega \right)$ are defined by
\begin{multline*}
D^{(1)}_t F = f^{(1)}(t) \\
+ \sum_{n = 1}^{\infty} \sum_{(i_1,\dots,i_n)} \sum_{k=1}^n \ind_{\{i_k = 1\}}  J_{(i_1,\dots,\widehat{i}_k,\dots,i_n)} \bigr( f^{(i_1,\dots,i_n)} (\underbrace{\dots}_{k-1},t,\underbrace{\dots}_{n - k}) \ind_{\Sigma^k_{n-1} (t)} \bigr)
\end{multline*}
and
\begin{multline*}
D^{(2)}_{t,z} F = f^{(2)}(t,z) \\
+ \sum_{n = 1}^{\infty} \sum_{(i_1,\dots,i_n)} \sum_{k=1}^n \ind_{\{i_k = 2\}}  J_{(i_1,\dots,\widehat{i}_k,\dots,i_n)} \bigr( f^{(i_1,\dots,i_n)} (\underbrace{\dots}_{k-1},(t,z),\underbrace{\dots}_{n - k}) \ind_{\Sigma^k_{n-1} (t)} \bigr)
\end{multline*}
if $F = \esp[F] + \sum_{n = 1}^{\infty} \sum_{(i_1,\dots,i_n)} J_{(i_1,\dots,i_n)} \left( f^{(i_1,\dots,i_n)} \right)$ is in $\mathbb{D}^{(1)}$ or $\mathbb{D}^{(2)}$.
\end{defn}

\begin{rem}\label{R:commutativity}
For an iterated integral, the Malliavin derivatives have a property similar to the classical commutativity relationship. Indeed, if $F = J_{(i_1,\dots,i_n)} (f)$, then
\begin{equation*}
D^{(2)}_{t,z} F = \int_t^T D^{(2)}_{t,z} J_{(i_1,\dots,i_{n-1})} \bigr( f(\cdot,s) \ind_{\Sigma_{n-1} (s)} \bigr) \, W(ds)
\end{equation*}
if $i_n = 1$ and
\begin{multline*}
D^{(2)}_{t,z} F = J_{(i_1,\dots,i_{n-1})} \bigr( f(\cdot,(t,z)) \ind_{\Sigma_{n-1} (t)} \bigr) \\
+ \int_t^T \int_{\reels} D^{(2)}_{t,z} J_{(i_1,\dots,i_{n-1})} \bigr( f(\cdot,(s,y)) \ind_{\Sigma_{n-1} (s)} \bigr) \, \Ntildedyds
\end{multline*}
if $i_n = 2$. A similar result holds for $D^{(1)} F$. 
\end{rem}

\begin{rem}\label{R:extension}
If $F = \esp[F] + \sum_{n = 1}^{\infty} J_n (f_n)$, where $J_n = J_{(1,\dots,1)}$ is the iterated Brownian stochastic integral of order $n$, then
\begin{align*}
D^{(1)}_t F &= f_1(t) + \sum_{n = 2}^{\infty} \sum_{k=1}^n J_{n-1} \bigr( f_n (\cdot,t,\cdot) \ind_{\Sigma^k_{n-1} (t)} \bigr) \\
&= f_1(t) + \sum_{n = 2}^{\infty} J_{n-1} ( f_n (\cdot,t) ) ,
\end{align*}
because $\sum_{k=1}^n \ind_{\Sigma^k_{n-1} (t)} = \ind_{[0,T]} (t)$. This is the classical Brownian Malliavin derivative of $F$. The same extension clearly holds for the pure-jump case if the $1$'s are replaced by $2$'s.
\end{rem}

The definitions of $\mathbb{D}^{(1)}$ and $\mathbb{D}^{(2)}$ come from the fact that we want the codomains of $D^{(1)}$ and $D^{(2)}$ to be $L^2 \left( [0,T] \times \Omega \right)$ and $L^2 \left( [0,T] \times \reels \times \Omega \right)$ respectively. We finally define a norm for $DF = (D^{(1)} F, D^{(2)} F)$ in the following way:
$$
\| DF \|^2 = \| D^{(1)} F \|^2_{L^2(\mesleb \times \p)} + \| D^{(2)} F \|^2_{L^2(\mesleb \times \meslevy \times \p)} .
$$
This is a norm on the product space $L^2(\mesleb \times \p) \times L^2(\mesleb \times \meslevy \times \p)$.

\subsection{Properties and interpretation of the Malliavin derivatives}

We begin this section with a result concerned with the \textit{continuity} of $D$. It is an extension of Lemma 1.2.3 in Nualart \cite{nualart1995}. The proof is given in Appendix~\ref{A:lemtech2}.

\begin{lem}\label{L:lemtech2}
If $F$ belongs to $L^2(\Omega)$, if $(F_k)_{k \geq 1}$ is a sequence of elements in $\D$ converging to $F$ in the $L^2(\Omega)$-norm and if $\sup_{k \geq 1} \| D F_k \| < \infty$, then $F$ belongs to $\D$ and $(D F_k)_{k \geq 1}$ converges weakly to $D F$ in $L^2(\mesleb \times \p) \times L^2(\mesleb \times \meslevy \times \p)$.
\end{lem}

There is a similar and stronger result stated in \cite{lokka2004} (Lemma $6$); however we are unable to fill a gap in its proof.

The choice for the definitions of the Malliavin derivative operators was made to extend the classical Brownian Malliavin derivative as well as the Poisson random measure Malliavin derivative in a wider sense than Remark~\ref{R:extension}. As mentioned in the introduction, the classical Brownian Malliavin derivative can be defined by chaos expansions and as a weak derivative. In Nualart and Vives \cite{nualartvives1990}, it is proven that for the Poisson process there is an equivalence between the Malliavin derivative defined with chaos decompositions and another one defined by \textit{adding a mass} with a translation operator. This last result was extended by L{\o}kka \cite{lokka2004} to Poisson random measures. But now we will follow an idea of Le{\'o}n et al. \cite{leonetal2002} to prove that our derivative operators are extensions of the classical ones. Their method relies on the commutativity relationships between stochastic derivatives and stochastic integrals and on quadratic covariation for semimartingales; consequently, it is easily adaptable to our more general context. The details are given in Appendix~\ref{A:interpretation}.

\begin{thm}\label{T:extension}
On $\mathbb{D}^{(1)}$ the operator $D^{(1)}$ coincides with the Brownian Malliavin derivative and on $\mathbb{D}^{(2)}$ the operator $D^{(2)}$ coincides with the Poisson random measure Malliavin derivative.
\end{thm}

Hence, if $F \in \mathbb{D}^{(1)}$, all the results about the classical Brownian Malliavin derivative, such as the chain rule for Lipschitz functions, can be applied to $D^{(1)} F$; see Nualart \cite{nualart1995} for details. But this is also true for the Poisson random measure Malliavin derivative. For example, an important result in L{\o}kka \cite{lokka2004} is that if $F = g(X_{t_1},\dots,X_{t_n}) \in \mathbb{D}^{(2)}$ and
$$
(t,z) \mapsto g \left( X_{t_1} + z \ind_{[0,t_1]}(t), \dots, X_{t_n} + z \ind_{[0,t_n]}(t) \right) - g \left( X_{t_1}, \dots, X_{t_n} \right)
$$
belongs to $L^2(\mesleb \times \meslevy \times \p)$, then
$$
D^{(2)}_{t,z} F = g \left( X_{t_1} + z \ind_{[0,t_1]}(t), \dots, X_{t_n} + z \ind_{[0,t_n]}(t) \right) - g \left( X_{t_1}, \dots, X_{t_n} \right) .
$$
This is the \textit{adding a mass} formula. Consequently, it also applies in the context of a square-integrable Lévy process.

\subsection{A Clark-Ocone formula}

We now state and prove a Clark-Ocone type formula. This formula gives explicitly the integrands in the martingale representation of Theorem~\ref{T:prp} for a Malliavin-differentiable Lévy functional. It is interesting to note that no particular property of the directional derivatives are needed.

\begin{thm}\label{T:clarkoconeformula}
If $F$ belongs to $\D$, then
\begin{equation*}\label{E:clarkOconeLevy}
F = \esp[F] + \int_0^T \esp \bigr[ D^{(1)}_t F \mid \tribu_t \bigr] \, W(dt) + \intPoissonT \esp \bigr[ D^{(2)}_{t,z} F \mid \tribu_t \bigr] \, \Ntildedt .
\end{equation*}
\end{thm}
\begin{proof}
Suppose that $F$ has a chaos expansion given by
$$
F = \esp[F] + \sum_{n = 1}^{\infty} \sum_{(i_1,\dots,i_n)} J_{(i_1,\dots,i_n)} \bigr( f^{(i_1,\dots,i_n)} \bigr) .
$$
If for example we consider the derivative operator $D^{(2)}$, then from Remark~\ref{R:crpimpliesmrp} we have to show that
\begin{multline}\label{E:toprove}
\esp \bigr[ D^{(2)}_{t,z} F \mid \tribu_t \bigr] \\
= f^{(2)}(t,z) + \sum_{n = 1}^{\infty} \sum_{(i_1,\dots,i_n)} J_{(i_1,\dots,i_n)} \bigr( f^{(i_1,\dots,i_n,2)} (\cdot,(t,z)) \ind_{\Sigma_n (t)} \bigr) .
\end{multline}
If $i_k = 2$, then
\begin{multline*}
\esp \bigr[ J_{(i_1,\dots,\widehat{i}_k,\dots,i_n)} \bigr( f^{(i_1,\dots,i_n)} (\cdot,(t,z),\cdot) \ind_{\Sigma^k_{n-1} (t)} \bigr) \bigr| \tribu_t \bigr] \\ =
\begin{cases}
0 & \text{if $k = 1,2,\dots,n-1$;}\\
J_{(i_1,\dots,i_{n-1})} \bigr( f^{(i_1,\dots,i_{n-1},2)} (\cdot,(t,z)) \ind_{\Sigma_{n-1} (t)} \bigr) & \text{if $k = n$,}
\end{cases}
\end{multline*}
because when $k = 1,2,\dots,n-1$ the outermost stochastic integral in the iterated integral $J_{(i_1,\dots,\widehat{i}_k,\dots,i_n)}$ starts after time $t$. By the definition of $ D^{(2)}_{t,z} F$, this implies that Equation~\eqref{E:toprove} is satisfied. The same argument works for the derivative operator $D^{(1)}$ and thus the result follows.
\end{proof}

\section{Martingale representation of the maximum}

Our main goal was to provide a detailed construction of a chaotic Malliavin derivative and a Clark-Ocone formula. Now, to illustrate the results, we compute the explicit martingale representation of the maximum of the Lévy process $X$.

For $0 \leq s < t \leq T$, define $M_{s,t} = \sup_{s \leq r \leq t} X_r$ and $M_t = M_{0,t}$. If $\esp [M_T] < \infty$, then one can show that
\begin{equation}\label{E:shiryaev}
\esp [M_T \mid \tribu_t] = M_t + \int_{M_t - X_t}^{\infty} \bar{F}_{T-t}(z) \, dz ,
\end{equation}
where $\bar{F}_{s}(z) = \p \{ M_s > z \}$; see Shiryaev and Yor \cite{shiryaevyor2004} and Graversen et al. \cite{graversenetal2001}. We will use this equality to prove the next proposition.

\begin{prop}
If $X$ a square-integrable Lévy process with Lévy-Itô decomposition
$$
X_t = \mu t + \sigma W_t + \intPoissont z \, \Ntildeds ,
$$
then its running maximum admits the following martingale representation:
$$
M_T = \esp [M_T] + \int_0^T \phi(t) \, W(dt) + \int_0^T \int_{\reels} \psi(t,z) \, \Ntildedt
$$
with $\phi(t) = \sigma \bar{F}_{T-t}(a)$ and $\psi(t,z) = \esp \left[ \left( M_{T-t} + z - a \right)^+ \right] - \int_a^{\infty} \bar{F}_{T-t}(x) \, dx$, where $a = M_t - X_t$.
\end{prop}
\begin{proof}
Since $X$ is a square-integrable martingale with drift, from Doob's maximal inequality we have that $M_T$ is a square-integrable random variable; see Theorem $20$ in Protter \cite{protter2004}. Let $(t_k)_{k \geq 1}$ be a dense subset of $[0,T]$, let $F = M_T$ and, for each $n \geq 1$, define $F_n = \max \{ X_{t_1}, \dots, X_{t_n} \}$. Clearly, $(F_n)_{n \geq 1}$ is an increasing sequence bounded by $F$. Hence, $F_n$ converges to $F$ in the $L^2(\Omega)$-norm when $n$ goes to infinity.

We want to prove that each $F_n$ is Malliavin differentiable, i.e. that each $F_n$ belongs to $\D = \mathbb{D}^{(1)} \cap \mathbb{D}^{(2)}$. This follows from the following two facts. First, since
$$
(x_1,\dots,x_n) \mapsto \max \{ x_1, \dots, x_n \}
$$
is a Lipschitz function on $\reels^n$ and since $D^{(1)}$ behaves like the classical Brownian Malliavin derivative on the Brownian part of $F_n$, we have that
$$
0 \leq D^{(1)}_t F_n = \sum_{k=1}^n \sigma \ind_{\{t \leq t_k\}} \ind_{A_k} \leq \sum_{k=1}^n \sigma \ind_{A_k} = \sigma ,
$$
where $A_1 = \{F_n = X_{t_1}\}$ and $A_k = \{F_n \neq X_{t_1} , \dots, F_n \neq X_{t_{k-1}}, F_n = X_{t_k}\}$ for $2 \leq k \leq n$. This implies that $\sup_{n \geq 1} \| D^{(1)} F_n \|_{L^2([0,T] \times \Omega)} \leq \sigma^2 T$. Secondly, since $D^{(2)}$ behaves like the Poisson random measure Malliavin derivative on the Poisson part of $F_n$, we have that
$$
0 \leq \bigr| D^{(2)}_{t,z} F_n \bigr| = \bigr| \max \left\lbrace X_{t_1} + z \ind_{\{t < t_1\}}, \dots, X_{t_n} + z \ind_{\{t < t_n\}} \right\rbrace - F_n \bigr| \leq |z| ,
$$
where the equality is justified by the following inequality:
\begin{multline*}
\bigr\| \max \left\lbrace X_{t_1} + z \ind_{\{t < t_1\}}, \dots, X_{t_n} + z \ind_{\{t < t_n\}} \right\rbrace - F_n \bigr\|_{L^2([0,T] \times \reels \times \Omega)}^2 \\
\leq T \int_{\reels} z^2 \, \meslevy(dz) .
\end{multline*}
Indeed, if $z \geq 0$, then
$$
0 \leq \max \left\lbrace X_{t_1} + z \ind_{\{t < t_1\}}, \dots, X_{t_n} + z \ind_{\{t < t_n\}} \right\rbrace - F_n \leq z ,
$$
and, if $z < 0$, then
\begin{align*}
0 & \leq F_n - \max \left\lbrace X_{t_1} + z \ind_{\{t < t_1\}}, \dots, X_{t_n} + z \ind_{\{t < t_n\}} \right\rbrace \\
& = F_n + \min \left\lbrace - X_{t_1} + |z| \ind_{\{t < t_1\}}, \dots, - X_{t_n} + |z| \ind_{\{t < t_n\}} \right\rbrace \\
& = \min \left\lbrace F_n - X_{t_1} + |z| \ind_{\{t < t_1\}}, \dots, F_n - X_{t_n} + |z| \ind_{\{t < t_n\}} \right\rbrace \\
& \leq |z| .
\end{align*}
This implies that $\sup_{n \geq 1} \| D^{(2)} F_n \|_{L^2([0,T] \times \reels \times \Omega)} \leq T \int_{\reels} z^2 \, \meslevy(dz)$.

Consequently, $\sup_{n \geq 1} \| D F_n \|^2 \leq T (\sigma^2 + \int_{\reels} z^2 \, \meslevy(dz))$ and by Theorem~\ref{T:extension} we have that $F$ is Malliavin differentiable. By the uniqueness of a weak limit, this means that taking the limit of $D^{(1)}_t F_n$ when $n$ goes to infinity yields
$$
D^{(1)}_t F = \sigma \ind_{[0, \tau]}(t) ,
$$
where $\tau$ is the first random time when the Lévy process $X$ (not the Brownian motion $W$) reaches its supremum on $[0,T]$, and
$$
D^{(2)}_{t,z} F = \sup_{0 \leq s \leq T} \left( X_s + z \ind_{\{t < s \}} \right) - M_T .
$$
Hence,
\begin{align*}
\esp \left[ D^{(1)}_t F \mid \tribu_t \right] &= \sigma \p \left\lbrace M_t < M_{t,T} \mid \tribu_t \right\rbrace\\
&= \sigma \p \left\lbrace M_{T-t} > a \right\rbrace ,
\end{align*}
where $a = M_t - X_t$. Since $M_{t,T} - X_t$ is independent of $\tribu_t$ and has the same law as $M_{T-t}$, then using Equation~\eqref{E:shiryaev} we get that
\begin{align*}
\esp \left[ D^{(2)}_{t,z} F \mid \tribu_t \right] &= \esp \left[ \sup_{0 \leq s \leq T} \left( X_s + z \ind_{\{t < s \}} \right) - M_T \mid \tribu_t \right] \\
&= \esp \left[ \max \{ M_t , M_{t,T} + z \} \mid \tribu_t \right] - \esp \left[ M_T \mid \tribu_t \right] \\
&= M_t + \esp \left[ \left( M_{t,T} + z - M_t \right)^+ \mid \tribu_t \right] - \esp \left[ M_T \mid \tribu_t \right] \\
&= \esp \left[ \left( M_{T-t} + z - a \right)^+ \right] - \int_a^{\infty} \bar{F}_{T-t}(x) \, dx .
\end{align*}
where $a = M_t - X_t$. The martingale representation follows from the Clark-Ocone formula of Theorem~\ref{T:clarkoconeformula}.
\end{proof}

This result extends the martingale representation of the running maximum of Brownian motion.

\section{Acknowledgements}

This paper is part of Jean-François Renaud's Ph.D. Thesis, written under the supervision of Bruno Rémillard. We would like to thank the thesis external referee Wim Schoutens for his careful reading and also Manuel Morales for bringing \cite{petrou2006} to our attention. We also thank Martin Goldstein for fruitful comments.

Partial funding in support of this work was provided by a doctoral scholarship of the Institut de Finance Math{\'e}matique de Montr{\'e}al (IFM2) and a scholarship of the Institut de Sciences Math{\'e}matiques (ISM), as well as by grants from the Natural Sciences and Engineering Research Council of Canada (NSERC) and the Fonds qu\'e\-b\'e\-cois de la re\-cher\-che sur la na\-tu\-re et les tech\-no\-lo\-gies (FQRNT).

\appendix

\section{Proof of Lemma~\ref{L:lemtech2}}\label{A:lemtech2}

We have that
$$
\sup_{k \geq 1} \| D^{(1)} F_k \|_{L^2([0,T] \times \Omega)} < \infty
$$
and
$$
\sup_{k \geq 1} \| D^{(2)} F_k \|_{L^2([0,T] \times \reels \times \Omega)} < \infty .
$$
Since $L^2([0,T] \times \Omega)$ and $L^2([0,T] \times \reels \times \Omega)$ are reflexive Hilbert spaces, there exist a subsequence $(k_j)_{j \geq 1}$, an element $\alpha$ in $L^2([0,T] \times \Omega)$ and an element $\beta$ in $L^2([0,T] \times \reels \times \Omega)$ such that $D^{(1)} F_{k_j}$ converges to $\alpha$ in the weak topology of $L^2([0,T] \times \Omega)$ and $D^{(2)} F_{k_j}$ converges to $\beta$ in the weak topology of $L^2([0,T] \times \reels \times \Omega)$. Consequently, for any $h \in L^2([0,T])$, $g \in L^2([0,T] \times \reels)$ and $f \in L^2(\sigman)$, we have that
$$
\left\langle D^{(1)} F_{k_j}, h \otimes J_{(i_1,\dots,i_n)}(f) \right\rangle_{L^2([0,T] \times \Omega)} \longrightarrow \left\langle \alpha, h \otimes J_{(i_1,\dots,i_n)}(f) \right\rangle_{L^2([0,T] \times \Omega)}
$$
and
$$
\left\langle D^{(2)} F_{k_j}, g \otimes J_{(i_1,\dots,i_n)}(f) \right\rangle_{L^2([0,T] \times \reels \times \Omega)} \longrightarrow \left\langle \beta, g \otimes J_{(i_1,\dots,i_n)}(f) \right\rangle_{L^2([0,T] \times \reels \times \Omega)}
$$
when $j$ goes to infinity.

Let $F = \esp[F] + \sum_{n = 1}^{\infty} \sum_{(i_1,\dots,i_n)} J_{(i_1,\dots,i_n)} ( f^{(i_1,\dots,i_n)} )$ and $F_{k_j} = \esp[F_{k_j}] + \sum_{n = 1}^{\infty} \sum_{(i_1,\dots,i_n)} J_{(i_1,\dots,i_n)} ( f_{k_j}^{(i_1,\dots,i_n)} )$ be the chaos representations of $F$ and $F_{k_j}$. By definition, we have that
\begin{multline}\label{E:deriveeun}
D^{(1)}_t F_{k_j} = f_{k_j}^{(1)}(t) \\
+ \sum_{n = 1}^{\infty} \sum_{(i_1,\dots,i_n)} \sum_{k=1}^n \ind_{\{i_k = 1\}}  J_{(i_1,\dots,\widehat{i}_k,\dots,i_n)} \bigr( f_{k_j}^{(i_1,\dots,i_n)} (\cdot,t,\cdot) \ind_{\Sigma^k_{n-1} (t)} \bigr) .
\end{multline}

By the linearity of the iterated integrals, the convergence of $F_{k_j}$ toward $F$ implies that
\begin{multline*}
\left\| \sum_{n = 1}^{\infty} \sum_{(i_1,\dots,i_n)} J_{(i_1,\dots,i_n)} \left( f^{(i_1,\dots,i_n)} - f_{k_j}^{(i_1,\dots,i_n)} \right) \right\|_{L^2(\Omega)} \\
= \sum_{n = 1}^{\infty} \sum_{(i_1,\dots,i_n)} \left\| f^{(i_1,\dots,i_n)} - f_{k_j}^{(i_1,\dots,i_n)} \right\|^2_{L^2(\sigman)}
\end{multline*}
goes to $0$ when $k$ tends to infinity. Consequently, it implies that each $f_{k_j}^{(i_1,\dots,i_n)}$ converges to $f^{(i_1,\dots,i_n)}$ when $j$ goes to infinity. So, using Proposition~\ref{P:orthogonality} and the expression of the derivative in Equation~\eqref{E:deriveeun}, we get that
\begin{align*}
& \bigr\langle D^{(1)} F_{k_j}, h \otimes J_{(i_1,\dots,i_n)}(f) \bigr\rangle_{L^2([0,T] \times \Omega)} \\
&= \sum_{k=1}^{n+1} \int_0^T \esp \bigr[ J_{(i_1,\dots,i_n)} \bigr( f_{k_j}^{(i_1,\dots,i_{k-1},1,i_k,\dots,i_n)} (\cdot,t,\cdot) \ind_{\Sigma^k_n (t)} \bigr) J_{(i_1,\dots,i_n)}(f) \bigr] h(t) \, dt \\
&= \sum_{k=1}^{n+1} \int_0^T \bigr\langle f_{k_j}^{(i_1,\dots,i_{k-1},1,i_k,\dots,i_n)} (\cdot,t,\cdot) \ind_{\Sigma^k_n (t)}, f \bigr\rangle_{L^2(\sigman)} h(t) \, dt
\end{align*}
and, as $j$ goes to infinity, this quantity tends to
$$
\sum_{k=1}^{n+1} \int_0^T \bigr\langle f^{(i_1,\dots,i_{k-1},1,i_k,\dots,i_n)} (\cdot,t,\cdot) \ind_{\Sigma^k_n (t)}, f \bigr\rangle_{L^2(\sigman)} h(t) \, dt .
$$
This holds for any multi-index $(i_1,\dots,i_n)$ and functions $h$ and $f$. Consequently,
\begin{multline*}
\alpha(t) = f^{(1)}(t) \\
+ \sum_{n = 1}^{\infty} \sum_{(i_1,\dots,i_n)} \sum_{k=1}^n \ind_{\{i_k = 1\}}  J_{(i_1,\dots,\widehat{i}_k,\dots,i_n)} \bigr( f^{(i_1,\dots,i_n)} (\cdot,t,\cdot) \ind_{\Sigma^k_{n-1} (t)} \bigr)
\end{multline*}
and $F$ belongs to $\mathbb{D}^{(1)}$ with $D^{(1)} F = \alpha$ by the unicity of the weak limit. Moreover, for any weakly convergent subsequence the limit must be equal to $D^{(1)} F$ and this implies the weak convergence of the whole sequence. The same argument works to prove that $F$ belongs to $\mathbb{D}^{(2)}$ and that $(D^{(2)} F_k)_{k \geq 1}$ converges weakly to $D^{(2)} F$ in $L^2(\mesleb \times \meslevy \times \p)$.

\section{Interpretation of the directional derivatives}\label{A:interpretation}

We consider the product probability space
$$
\left( \Omega_W \times \Omega_N, \tribu_W \times \tribu_N, \p_W \times \p_N \right)
$$
which is the product of the canonical space of the Brownian motion $W$ and the canonical space of the pure-jump Lévy process
$$
N_t = \int_0^t \int_{\reels} z \, \Ntildeds
$$
previously defined in Equation~\eqref{E:defN}; see Solé et al. \cite{soleetal2007} for more details on this last canonical space. Since $L^2(\Omega_W \times \Omega_N)$ is isometric to $L^2(\Omega_W ; L^2(\Omega_N))$ and to $L^2(\Omega_N ; L^2(\Omega_W))$ as Hilbert spaces, we will use the theory of the Brownian Malliavin derivative and the Poisson random measure Malliavin derivative for Hilbert-valued random variables (see \cite{nualart1995} and \cite{nualartvives1990}). This is possible because both operators are closable.

The Brownian Malliavin derivative for Hilbert-valued random variables will be denoted by $D^W$ and the Poisson random measure Malliavin derivative for Hilbert-valued random variables by $D^N$. If we define $\widetilde{W} = (\widetilde{W}_t)_{t \in [0,T]}$ on $\Omega_W \times \Omega_N$ by
$$
\widetilde{W}_t (\omega,\omega^{\prime}) = \omega(t)
$$
and $\widetilde{N} = (\widetilde{N}_t)_{t \in [0,T]}$ by
$$
\widetilde{N}_t (\omega,\omega^{\prime}) = \omega^{\prime}(t) ,
$$
then the process $\widetilde{X}_t = \mu t + \sigma \widetilde{W}_t + \widetilde{N}_t$ has the same distribution as our initial Lévy process $X_t = \mu t + \sigma W_t + N_t$. For notational simplicity, in what follows we will write $W_t(\omega)$ and $N_t(\omega^{\prime})$ instead of $\widetilde{W}_t (\omega,\omega^{\prime})$ and $\widetilde{N}_t (\omega,\omega^{\prime})$ respectively.

We will proceed by induction. If $F = \int_0^T f(t) \, W(dt)$, then clearly
\begin{equation*}
D^{(1)}_t F = D^W_t F = f(t) \quad \text{and} \quad D^{(2)}_{t,z} F = D^N_{t,z} F = 0 ,
\end{equation*}
while if $G = \int_0^T \int_{\reels} g(t,z) \, \Ntildedt$, then
\begin{equation*}
D^{(1)}_t G = D^W_t G = 0 \quad \text{and} \quad D^{(2)}_{t,z} G = D^N_{t,z} G = g(t,z) .
\end{equation*}

Thus, for a fixed $n \geq 1$, we assume that $D^{(1)}$ and $D^W$ coincide for any random variable with chaos expansion of order $n$. First, let $F$ be of the form
$$
F = J_{(i_1,\dots,i_n,1)} (f_1 \otimes \dots \otimes f_n \otimes f_{n+1})  = \int_0^T g(s) f_{n+1}(s) \, W(ds) ,
$$
where
\begin{align}\label{E:g}
g(s) =  J_{(i_1,\dots,i_n)} \bigr( f_1 \otimes \dots \otimes f_n \, \ind_{\Sigma_{n} (s)} \bigr) .
\end{align}
To ease the notation, $J_{(i_1,\dots,i_n)} ( f_1 \dots f_n )$ will mean $J_{(i_1,\dots,i_n)} ( f_1 \otimes \dots \otimes f_n )$. Using the commutativity relationship of Remark~\ref{R:commutativity} and the hypothesis of induction, we have that
\begin{align*}
D^{(1)}_t F &= f_{n+1}(t) g(t) + \int_t^T f_{n+1}(s) D^{(1)}_t g(s) \, W(ds) \\
&= f_{n+1}(t) g(t) + \int_t^T f_{n+1}(s) D^W_t g(s) \, W(ds) ,
\end{align*}
which is exactly $D^W_t F$, by the classical commutativity relationship of Equation~\eqref{E:commute}.

Secondly, now let $F$ be of the form
$$
F = J_{(i_1,\dots,i_n,2)} (f_1 \otimes \dots \otimes f_n \otimes f_{n+1})  = \int_0^T \int_{\reels} g(s-) f_{n+1}(s,z) \, \Ntildeds .
$$
We will use of the \textit{integration by parts formula} for semimartingales, that is
$$
\bigr[ Y^{(1)},Y^{(2)} \bigr]_t = Y^{(1)}_t Y^{(2)}_t - \int_0^t Y^{(1)}_{s-} \, dY^{(2)}_s - \int_0^t Y^{(2)}_{s-} \, dY^{(1)}_s
$$
if $Y^{(1)}$ and $Y^{(2)}$ are semimartingales; see Protter \cite{protter2004} for details. If $Y^{(1)}_t = g(t)$ and $Y^{(2)}_t = \int_0^t \int_{\reels} f_{n+1}(s,z) \, \Ntildeds$, we get that
\begin{multline*}
F = g(T) \int_0^T \int_{\reels} f_{n+1}(s,z) \, \Ntildeds \\
- \int_0^T \int_0^{t-} \int_{\reels} f_{n+1}(s,z) \, \Ntildeds dg(t) \\
- \left[ g(\cdot) , \int_0^{\cdot} \int_{\reels} f_{n+1}(s,z) \, \Ntildeds \right]_T .
\end{multline*}

We now consider the two cases where $i_n = 1$ and $i_n = 2$ separately. We have that
\begin{equation*}
g(t) =
\begin{cases}
\int_0^t h(s) f_n(s) \, W(ds) & \text{if $i_n = 1$;} \\
\int_0^t \int_{\reels} h(s-) f_n(s,z) \, \Ntildeds & \text{if $i_n = 2$,}
\end{cases}
\end{equation*}
where $h(s) =  J_{(i_1,\dots,i_n)} \bigr( f_1 \otimes \dots \otimes f_{n-1} \, \ind_{\Sigma_{n-1} (s)} \bigr)$. If $i_n = 1$,  then
\begin{multline*}
F = g(T) \int_0^T \int_{\reels} f_{n+1}(t,z) \, \Ntildedt \\
- \int_0^T \left[ \int_0^t \int_{\reels} f_{n+1}(s,y) \, \Ntildedyds \right] h(t) f_n(t) \, W(dt) .
\end{multline*}
If $i_n = 2$, then
\begin{multline*}
F = g(T) \int_0^T \int_{\reels} f_{n+1}(t,z) \, \Ntildedt \\
- \int_0^T \int_{\reels} \left[ \int_0^{t-} \int_{\reels} f_{n+1}(s,y) \, \Ntildedyds \right] h(t-) f_n(t,z) \, \Ntildedt \\
- \int_0^T \int_{\reels} h(t-) f_n(t,z) f_{n+1}(t,z) \, N(dt,dz) .
\end{multline*}
Note that the last term is an iterated integral of order $n$ (with respect to $N(dt,dz)$ for the outermost integral, not $\Ntildedt$) since $h$ is an iterated integral of order $n-1$. So, by the hypothesis of induction, $D^{(1)}$ and $D^W$ agree for this functional. This is also true for $g(T)$. 

Consequently, we repeat the previous steps backward with $D^{(1)}$. If $i_n = 1$, then
\begin{align*}
D^W_t F &= \left( D^W_t g(T) \right)  \int_0^T \int_{\reels} f_{n+1}(s,y) \, \Ntildedyds \\
& \qquad - h(t) f_n(t) \int_0^t \int_{\reels} f_{n+1}(s,y) \, \Ntildedyds \\
& \qquad \qquad - \int_t^T \left[ \int_0^s \int_{\reels} f_{n+1}(r,y) \, \Ntildedydr \right] \left( D^W_t h(s) \right) f_n(s) \, W(ds) \\
&= \left( D^{(1)}_t g(T) \right)  \int_0^T \int_{\reels} f_{n+1}(s,y) \, \Ntildedyds \\
& \qquad - h(t) f_n(t) \int_0^t \int_{\reels} f_{n+1}(s,y) \, \Ntildedyds \\
& \qquad \qquad - \int_t^T \left[ \int_0^s \int_{\reels} f_{n+1}(r,y) \, \Ntildedydr \right] \left( D^{(1)}_t h(s) \right) f_n(s) \, W(ds) \\
&= D^{(1)}_t \left( g(T) \int_0^T \int_{\reels} f_{n+1}(s,y) \, \Ntildedyds \right) \\
& \qquad - D^{(1)}_t \int_0^T \left[ \int_0^s \int_{\reels} f_{n+1}(r,y) \, \Ntildedydr \right] h(s) f_n(s) \, W(ds) \\
&= D^{(1)}_t \biggr( g(T) \int_0^T \int_{\reels} f_{n+1}(s,y) \, \Ntildedyds \\
& \qquad \qquad \qquad - \int_0^T \left[ \int_0^s \int_{\reels} f_{n+1}(r,y) \, \Ntildedydr \right] \, dg(s) \biggr) \\
&= D^{(1)}_t F ,
\end{align*}
and if $i_n = 2$, then the same steps are valid since $D^W$ and $D^{(1)}$ coincide on the extra term.

The equivalence between $D^{(1)}$ and $D^W$ follows from the following fact: for a fixed $n \geq 1$ and a fixed multi-index $(i_1,\dots,i_n)$, the linear subspace of $L^2(\sigman)$ generated by functions of the form
\begin{equation}\label{E:tensor}
f_1 \otimes \dots \otimes f_n ,
\end{equation}
is dense. Indeed, for $f \in L^2(\sigman)$, there exists a sequence $(f_n)_{n \geq 1}$, whose elements are finite sums of functions as in Equation~\eqref{E:tensor}, that converges to $f$. We know that $D^{(1)}$ and $D^W$ are equal for each $f_n$. Since $D^{(1)}$ and $D^W$ are continuous (see Lemma~\ref{L:lemtech2}), they also coincide for $f$.

We can apply the same machinery to show that $D^N$ and $D^{(2)}$ are the same.

\bibliographystyle{abbrv}
\bibliography{renaudremillard_arxiv}

\begin{thebibliography}{10}

\bibitem{aaseetal2000}
K.~Aase, B.~{\O}ksendal, N.~Privault, and J.~Ub{\o}e.
\newblock White noise generalizations of the {C}lark-{H}aussmann-{O}cone
  theorem with application to mathematical finance.
\newblock {\em Finance Stoch.}, 4(4):465--496, 2000.

\bibitem{ballyetal2007}
V.~Bally, M.-P. Bavouzet, and M.~Messaoud.
\newblock Integration by parts formula for locally smooth laws and applications
  to sensitivity computations.
\newblock {\em Ann. Appl. Probab.}, 17(1):33--66, 2007.

\bibitem{bavouzetmessaoud2006}
M.-P. Bavouzet-Morel and M.~Messaoud.
\newblock Computation of {G}reeks using {M}alliavin's calculus in jump type
  market models.
\newblock {\em Electron. J. Probab.}, 11, 2006.

\bibitem{benthetal2003}
F.~E. Benth, G.~Di~Nunno, A.~L{\o}kka, B.~{\O}ksendal, and F.~Proske.
\newblock Explicit representation of the minimal variance portfolio in markets
  driven by {L}\'evy processes.
\newblock {\em Math. Finance}, 13(1):55--72, 2003.

\bibitem{bertoin1996}
J.~Bertoin.
\newblock {\em L\'evy processes}.
\newblock Cambridge University Press, 1996.

\bibitem{davisjohansson2006}
M.~H.~A. Davis and M.~P. Johansson.
\newblock Malliavin {M}onte {C}arlo {G}reeks for jump diffusions.
\newblock {\em Stochastic Process. Appl.}, 116(1):101--129, 2006.

\bibitem{dermoune1990}
A.~Dermoune.
\newblock Distributions sur l'espace de {P}. {L}\'evy et calcul stochastique.
\newblock {\em Ann. Inst. H. Poincar\'e Probab. Statist.}, 26(1):101--119,
  1990.

\bibitem{graversenetal2001}
S.~E. Graversen, G.~Peskir, and A.~N. Shiryaev.
\newblock Stopping {B}rownian motion without anticipation as close as possible
  to its ultimate maximum.
\newblock {\em Theory of Probability and its Applications}, 45(1):125--136,
  2001.

\bibitem{ito1956}
K.~It{\^o}.
\newblock Spectral type of the shift transformation of differential processes
  with stationary increments.
\newblock {\em Trans. Amer. Math. Soc.}, 81:253--263, 1956.

\bibitem{kulik2006}
A.~M. Kulik.
\newblock Malliavin calculus for {L}{\'e}vy processes with arbitrary {L}{\'e}vy
  measures.
\newblock {\em Theory of Probability and its Applications}, (72):75--92, 2006.

\bibitem{kunita2004}
H.~Kunita.
\newblock Representation of martingales with jumps and applications to
  mathematical finance.
\newblock In {\em Stochastic analysis and related topics in Kyoto}, pages
  209--232. Math. Soc. Japan, 2004.

\bibitem{kunitawatanabe1967}
H.~Kunita and S.~Watanabe.
\newblock On square integrable martingales.
\newblock {\em Nagoya Math. J.}, 30:209--245, 1967.

\bibitem{leonetal2002}
J.~A. Le{\'o}n, J.~L. Sol{\'e}, F.~Utzet, and J.~Vives.
\newblock On {L}\'evy processes, {M}alliavin calculus and market models with
  jumps.
\newblock {\em Finance Stoch.}, 6(2):197--225, 2002.

\bibitem{lokka2004}
A.~L{\o}kka.
\newblock Martingale representation of functionals of {L}\'evy processes.
\newblock {\em Stochastic Anal. Appl.}, 22(4):867--892, 2004.

\bibitem{lytvynov2003}
E.~Lytvynov.
\newblock Orthogonal decompositions for {L}\'evy processes with an application
  to the gamma, {P}ascal, and {M}eixner processes.
\newblock {\em Infin. Dimens. Anal. Quantum Probab. Relat. Top.}, 6(1):73--102,
  2003.

\bibitem{maetal1998}
J.~Ma, P.~Protter, and J.~San~Martin.
\newblock Anticipating integrals for a class of martingales.
\newblock {\em Bernoulli}, 4(1):81--114, 1998.

\bibitem{meyer1976}
P.-A. Meyer.
\newblock Un cours sur les int\'egrales stochastiques.
\newblock In {\em S\'eminaire de Probabilit\'es, X}, volume 511 of {\em Lecture
  Notes in Math.}, pages 245--400. Springer, 1976.

\bibitem{nualart1995}
D.~Nualart.
\newblock {\em The {M}alliavin calculus and related topics}.
\newblock Springer-Verlag, 1995.

\bibitem{nualartschoutens2000}
D.~Nualart and W.~Schoutens.
\newblock Chaotic and predictable representations for {L}\'evy processes.
\newblock {\em Stochastic Process. Appl.}, 90(1):109--122, 2000.

\bibitem{nualartvives1990}
D.~Nualart and J.~Vives.
\newblock Anticipative calculus for the {P}oisson process based on the {F}ock
  space.
\newblock In {\em S\'eminaire de Probabilit\'es, XXIV, 1988/89}, volume 1426 of
  {\em Lecture Notes in Math.}, pages 154--165. Springer, 1990.

\bibitem{petrou2006}
E.~Petrou.
\newblock Malliavin calculus in {L}\'evy spaces and applications in finance.
\newblock {\em preprint}, 2006.

\bibitem{protter2004}
P.~E. Protter.
\newblock {\em Stochastic integration and differential equations}.
\newblock Springer-Verlag, second edition, 2004.

\bibitem{schoutens2003}
W.~Schoutens.
\newblock {\em L\'evy processes in finance: pricing financial derivatives}.
\newblock Wiley, 2003.

\bibitem{shiryaevyor2004}
A.~N. Shiryaev and M.~Yor.
\newblock On stochastic integral representations of functionals of {B}rownian
  motion.
\newblock {\em Theory of Probability and its Applications}, 48(2):304--313,
  2004.

\bibitem{soleetal2007}
J.~L. Sol{\'e}, F.~Utzet, and J.~Vives.
\newblock Canonical {L}\'evy process and {M}alliavin calculus.
\newblock {\em Stochastic Process. Appl.}, 117(2):165--187, 2007.

\end{thebibliography}

\end{document}